\documentclass[11pt]{article}
\usepackage{amsfonts}
\usepackage[fleqn]{amsmath}
\usepackage{amssymb}
\usepackage{dsfont}
\usepackage[pdftex=true,colorlinks=true,plainpages=false]{hyperref}
    \hypersetup{
    colorlinks=true,%
    citecolor=black,%
    filecolor=black,%
    linkcolor=black,%
    urlcolor=black,%
    }
\usepackage[top=1.75cm,bottom=1.75cm,left=2cm,right=2cm]{geometry}
\usepackage{setspace}
\usepackage{enumitem}
\usepackage{float}
    \restylefloat{table}
\usepackage{graphicx}
\usepackage{epstopdf}
\usepackage{caption}
\usepackage{subcaption}

\renewenvironment{abstract}
{\small
\begin{center}
\bfseries\abstractname\vspace{-.5em}\vspace{0pt}
\end{center}
\list{}{%
\setlength{\leftmargin}{0.75cm}%
\setlength{\rightmargin}{\leftmargin}%
}%
\item\relax}
{\endlist}

\newtheorem{thm}{Theorem}[section]
\newtheorem{lem}{Lemma}[section]
\newtheorem{prop}{Proposition}[section]
\newtheorem{rem}{Remark}[section]

\numberwithin{equation}{section}

\renewcommand{\l}{\linebreak}

\title{Rapid stabilization for a wave equation with boundary disturbance}

\author{Patricio Guzm\'an\footnotemark[1],\hspace{0.275cm}Agust\'in Huerta\footnotemark[1]\hspace{0.275cm}and Hugo Parada\footnotemark[2]}

\begin{document}
\setlength{\abovedisplayskip}{0pt}

\footnotetext[1]{Departamento de Matem\'atica, Universidad T\'ecnica Federico Santa Mar\'ia, Valpara\'iso, Chile.\\
E-mail: {\tt patricio.guzmanm@usm.cl, huertasotoagustin@gmail.com (corresponding author)}}

\footnotetext[2]{Universit\'e de Lorraine, CNRS, Inria, IECL, F-54000 Nancy, France.\\
E-mail: {\tt hugo.parada@inria.fr}}

\date{}

\maketitle

\begin{abstract}
In this paper we study the rapid stabilization of an unstable wave equation, in which an unknown\l disturbance is located at the boundary condition. We address two different boundary conditions: Dirichlet-Dirichlet and Dirichlet-Neumann. In both cases, we design a feedback law, located at the same place as the unknown disturbance, that forces the exponential decay of the energy for any desired decay rate while suppressing the effects of the unknown disturbance. For the feedback design we employ the backstepping method, Lyapunov techniques and the sign multivalued operator. The well-posedness of the closed-loop system, which is a differential inclusion, is shown with the maximal monotone operator theory.
\end{abstract}

\vspace{0.3cm}


\vspace{0.5cm}



\section{Introduction}

\hspace{0.75cm}\textbf{Some general comments.} Partial differential equations have been widely and successfully used to derive mathematical models for diverse phenomena, such as the temperature on solids or fluids \cite{cg2015,k2018,s2022}, the propagation of waves in a medium \cite{dv2013,g2016,w2019}, the deflection of strings or beams \cite{ll1988,r2019}, to name a few remarkable examples. Once a mathematical model is established one relevant task in control theory is to design feedback laws to stabilize the state of the system to their equilibria or to another state of interest \cite{cz1995,lgm1999,k2002,l2010}. In the literature we can frequently find that mathematical models are analyzed under ideal assumptions in which disturbances are neglected for the sake of simplicity, but disturbances are always present and indeed may correspond to an additional source of instability. Accordingly, it is relevant to include disturbances in the study of the stabilization problem.\\

\textbf{A wave equation with boundary disturbance.} Let $L\in(0,\infty)$. Let $\lambda$, $\alpha$ and $\beta$ be time\l independent coefficients. Let $p=p(t)$ represent an unknown disturbance and let $U=U(t)$ represent a control to be designed. Let $(u_{0},u_{1})$ be an initial condition. Then, as in \cite{sck2010}, let us consider the following two cases.

\newpage

$\bullet$ Dirichlet-Dirichlet case:

\begin{equation}
\label{wave_dd}
\left\{\begin{array}{rl}
u_{tt}=u_{xx}+2\lambda(x)u_{t}+\alpha(x)u_{x}+\beta(x)u,&~(t,x)\in(0,\infty)\times(0,L),\\[1ex]
u(t,0)=0,~u(t,L)=p(t)+U(t),&~t\in(0,\infty),\\[1ex]
u(0,x)=u_{0}(x),~u_{t}(0,x)=u_{1}(x),&~x\in(0,L).
\end{array}\right.
\end{equation}

$\bullet$ Dirichlet-Neumann case:

\begin{equation}
\label{wave_dn}
\left\{\begin{array}{rl}
u_{tt}=u_{xx}+2\lambda(x)u_{t}+\alpha(x)u_{x}+\beta(x)u,&~(t,x)\in(0,\infty)\times(0,L),\\[1ex]
u(t,0)=0,~u_{x}(t,L)=p(t)+U(t),&~t\in(0,\infty),\\[1ex]
u(0,x)=u_{0}(x),~u_{t}(0,x)=u_{1}(x),&~x\in(0,L).
\end{array}\right.
\end{equation}

\textbf{Rapid stabilization.} For a regular enough state $(u,u_{t})$ of \eqref{wave_dd} or \eqref{wave_dn}, we define its energy by

\begin{equation}
\label{energy}
E(t)=\frac{1}{2}\int_{0}^{L}(|u_{t}|^{2}+|u_{x}|^{2})\,dx,~t\in[0,\infty).
\end{equation}

\noindent The purpose of this paper is to exponentially stabilize the regular enough state $(u,u_{t})$ to the rest position $(0,0)$ by means of a control $U$ that suppresses the effects of an unknown disturbance $p$. Being precise, we aim to prove the rapid stabilization, which in this paper is understood as follows:

\begin{equation}
\label{rapid}
\left\{\begin{array}{l}
\mbox{For any desired decay rate }d\in(0,\infty)\\[1ex]
\mbox{there exists a control }U(t)\\[1ex]
\mbox{and there exists a constant }C\in[1,\infty)\\[1ex]
\mbox{such that }E(t)\leq Ce^{-2d t}E(0)\mbox{ for all }t\in[0,\infty).
\end{array}\right.
\end{equation}

In virtue of \cite{l2003,sgk2009,sck2010,cc2013,cl2014,cl2015,rbcps2018,kg2019,glm2021,ccg2023} and also of many other works, as can be consulted in \cite{ks2008,vabk2026} for instance, we have that the backstepping method has been successfully applied to prove the rapid stabilization of several partial differential equations in the undisturbed case, which in this paper is when $p=0$.

\begin{rem}
In \cite{sck2010} it has been shown \eqref{rapid} for \eqref{wave_dd} and for \eqref{wave_dn} when $p=0$. In this paper we address the disturbed case. Then, if we denote by $U_{*}(t)$ the control designed in \cite{sck2010}, unfortunately we cannot ask $U(t)=U_{*}(t)-p(t)$ to solve the disturbed case, since we assume that $p(t)$ is unknown.
\end{rem}

In this paper we employ the backstepping transformation constructed in \cite{sck2010} and apply the idea\l behind the feedback design in \cite{gh2025} to show \eqref{rapid} for \eqref{wave_dd} and for \eqref{wave_dn}. Accordingly, although the\l disturbance is assumed unknown, we ask it to satisfy some assumptions, that shall be introduced next.\\

\textbf{Stabilization in the disturbed case.} The disturbed case has been object of recent interest, as can be consulted in Table \ref{literature}. There we present some of the concerned literature in the stabilization of partial differential equations subjected to unknown disturbances, acting either in the domain or at the boundary.

\newpage
\mbox{}
\vspace{-1.25cm}

\begin{table}[H]
\centering
\begin{tabular}{c|c|c}
Equation & Distributed disturbance & Boundary disturbance\\
\hline
Heat & \cite{gp2020} & \cite{kg2014,lw2015,gh2023,gh2025} \\
Wave & \cite{fx2015} & \cite{gk2012,gj2013,mgk2019} \\
Beam & \cite{a2023} & \cite{gk2014,jg2015,lx2017} \\
Schr\"odinger & -- & \cite{gl2014,kg2019}
\end{tabular}
\caption{\label{literature}Stabilization of partial differential equations subjected to unknown disturbances. Assumptions \textbf{(A1)} and \textbf{(A2)} for the unknown disturbance considered in this paper are based on the assumptions asked in the referenced works.}
\label{tl}
\end{table}

Given the works in Table \ref{literature}, it follows that the commonly accepted assumptions for the unknown disturbance involves knowing a uniform-in-time upper bound for its norm and also knowing its regularity. In this paper these assumptions are the following:\\

\noindent\textbf{(A1)} There exists a constant $P\in(0,\infty)$ such that $|p(t)|\leq P$ for every $t\in[0,\infty)$.

\noindent\textbf{(A2)$\mathbf{_{DD}}$} $p\in W^{2,1}(0,\infty)$ and $p(0)=0$.

\noindent\textbf{(A2)$\mathbf{_{DN}}$} $p\in W^{3,1}(0,\infty)$ and $p(0)=p'(0)=0$.\\

\noindent Some comments on these assumptions are the following.\\

The works in Table \ref{literature} use the knowledge of the uniform-in-time upper bound for the norm of the unknown disturbance for the design of controls as feedback laws. In that design it is commonly used the sign multivalued operator $\text{sign}:\mathbb{R}\to 2^{\mathbb{R}}$ ($2^{\mathbb{R}}$ denotes the power set of $\mathbb{R}$) given by

\begin{equation}
\label{sign}
\text{sign}(f)=\left\{\begin{array}{ccc}
\dfrac{f}{|f|} & \mbox{if} & f\neq0,\\[3ex]
[-1,1] & \mbox{if} & f=0,
\end{array}\right.
\end{equation}

\noindent to suppress the instability that may produce the unknown disturbance. To that end it is adequately used the property $|f|=\text{sign}(f)f$. These feedback laws also contain other (univalued) terms, but are not involved in the suppression.\\

The works in Table \ref{literature} use the knowledge of the regularity of the unknown disturbance to satisfy the required hypotheses to apply the well-posedness results for the corresponding closed-loop systems, which are differential inclusions due to \eqref{sign}. Regarding the well-posedness, only \cite{kg2014,gk2014,jg2015,kg2019} consider the notion of Filippov solution, requiring the introduction of a sliding surface, while the rest of the works consider the notion of classical solution, as understood in the maximal monotone operator theory. In this paper we consider the latter notion, beign the choice only due to mathematical convenience.\\

\textbf{Contribution of this paper.} Except by \cite{kg2019} (Schr\"odinger equation) and by \cite{gp2020,gh2025} (Heat equation),\l all the works in Table \ref{literature} prove the exponential decay of some energy for a certain decay rate that results from the proof. Accordingly, the obtained decay rate in those works is unchosen. Therefore, most of the results are not rapid in the sense described in \eqref{rapid}, being the case addressed in this paper (Wave equation) an open problem. To fill that gap, in this paper we employ the backstepping transformation constructed in \cite{sck2010} and apply the idea behind the feedback design in \cite{gh2025} to show the rapid stabilization of an unstable wave equation with boundary disturbance, meaning \eqref{rapid} for \eqref{wave_dd} and for \eqref{wave_dn}. The well-posedness results for the corresponding closed-loop systems are obtained thanks to the maximal monotone operator theory by adapting an optimization argument in \cite{cp1990}.

For the feedback design (see \cite{k1994,lgm1999,c2007,ks2008,l2010,m2013,gw2019} for instance) we take inspiration from \cite{gh2025}. Indeed, in \eqref{wave_dd} or \eqref{wave_dn} we consider $U(t)=U_{1}(t)+U_{2}(t)$, where $U_{1}(t)$ is used to achieve the desired decay\l rate and $U_{2}(t)$ is used to reject the effects of the unknown disturbance. Being precise, it is in the choice of $U_{1}(t)$ where the backstepping transformation constructed in \cite{sck2010} is employed and it is in the design of $U_{2}(t)$ where Lyapunov techniques and \eqref{sign} are adequately applied. Since the final expression for the feedback law $U(t)$ is multivalued because of \eqref{sign}, the corresponding closed-loop system is a differential inclusion, with the inclusion being located at the boundary condition. We prove the well-posedness result with the aid of the maximal monotone operator theory (see \cite{b1973,s1997,b2010} and/or \cite{bp1970,cp1972,ds1973,ps2010} for instance), in which we adapt an optimization argument in \cite{cp1990} that uses the Moreau regularization to treat the terms at the boundary condition. Our main results are given by Theorem \ref{dd_mr} (Dirichlet-Dirichlet case) and by Theorem \ref{dn_mr} (Dirichlet-Neumann case).\\

\textbf{Organization of this paper.} The rest of this paper is organized as follows. The feedback design for the Dirichlet-Dirichlet case (respectively: Dirichlet-Neumann case) is done in Section \ref{dd_fd} (respectively: Section \ref{dn_fd}), while the well-posedness of the corresponding closed-loop system is shown in Section \ref{dd_wp} (respectively: Section \ref{dn_wp}). The main result gathering the conclusions of both the feedback design and the well-posedness parts is stated in Section \ref{dd_c} (respectively: Section \ref{dn_c}). Finally, in Section \ref{cr} we provide closing comments.

\section{Dirichlet-Dirichlet case}

\hspace{0.75cm}In this section we present all the required analysis to prove Theorem \ref{dd_mr}, which is one of our main results and corresponds to the obtention of \eqref{rapid} for \eqref{wave_dd}. We start with a comment.

\begin{rem}
\label{simplification_dd}
As stated in \cite[Section 2]{sck2010}, we can transform \eqref{wave_dd} into another wave equation without\l first-order spatial derivative term by applying the transformation $v(t,x)=(e^{\frac{1}{2}\int_{0}^{x}\alpha(\tau)\,d\tau})u(t,x)$, thus\l obtaining

\begin{equation}
\label{dd_change}
\left\{\begin{array}{rl}
v_{tt}=v_{xx}+2\lambda(x)v_{t}+\left[-\frac{1}{2}\alpha'(x)-\frac{1}{4}\alpha(x)^{2}+\beta(x)\right]v,&~(t,x)\in(0,\infty)\times(0,L),\\[1ex]
v(t,0)=0,~v(t,L)=(e^{\frac{1}{2}\int_{0}^{L}\alpha(\tau)\,d\tau})p(t)+(e^{\frac{1}{2}\int_{0}^{L}\alpha(\tau)\,d\tau})U(t),&~t\in(0,\infty),\\[1ex]
v(0,x)=(e^{\frac{1}{2}\int_{0}^{x}\alpha(\tau)\,d\tau})u_{0}(x),~v_{t}(0,x)=(e^{\frac{1}{2}\int_{0}^{x}\alpha(\tau)\,d\tau})u_{1}(x),&~x\in(0,L).
\end{array}\right.
\end{equation}
\end{rem}

Thus, in this section we consider

\begin{equation}
\label{dd_1}
\left\{\begin{array}{rl}
u_{tt}=u_{xx}+2\lambda(x)u_{t}+\beta(x)u,&~(t,x)\in(0,\infty)\times(0,L),\\[1ex]
u(t,0)=0,~u(t,L)=p(t)+U(t),&~t\in(0,\infty),\\[1ex]
u(0,x)=u_{0}(x),~u_{t}(0,x)=u_{1}(x),&~x\in(0,L).
\end{array}\right.
\end{equation}

\subsection{Feedback design}
\label{dd_fd}

\hspace{0.75cm}In this section we prove the feedback design part of Theorem \ref{dd_mr}. To that end, we start by introducing all the required elements of the backstepping transformation constructed in \cite[Section 3]{sck2010}.\\

Let $d\in(0,\infty)$ be the desired decay rate of the exponential decay of the energy \eqref{energy}. We proceed to introduce some elements from \cite[Pages 4017, 4018 and 4019]{sck2010}, with $c(x)=d^{2}$ and $d(x)=d$ in their notation. For $x\in[0,L]$ let us consider

\begin{equation}
\label{BS_1}
\left\{\begin{array}{l}
a(x)=\lambda(x)+d,\\[1ex]
h(x)=\mbox{cosh}\left(\int_{0}^{x}a(\tau)\,d\tau\right),\\[1ex]
m(x)=\frac{h'(x)}{2a(x)}\left(2\lambda(x)+a(x)+a(0)\right)+\frac{h(x)}{2}\int_{0}^{x}\left(d^{2}-\lambda(y)^{2}-\beta(y)-d^{2}\right)\,dy.
\end{array}\right.
\end{equation}

\noindent Let $\mathcal{T}=\{(x,y)\in\mathbb{R}^{2}~/~0\leq y\leq x\leq L\}$. For $(x,y)\in\mathcal{T}$ let us consider

\begin{equation}
\label{BS_2}
\left\{\begin{array}{l}
\rho_{1}(x,y)=2\lambda(y)+2d,\\[1ex]
\rho_{2}(x,y)=d^{2}+\beta(y),\\[1ex]
\rho_{3}(x,y)=2\lambda(y)\beta(y)+2\lambda''(y)+2d\beta(y),\\[1ex]
\rho_{4}(x,y)=4\lambda'(y),\\[1ex]
\rho_{5}(x,y)=4\lambda(y)^{2}+4d\lambda(y)+d^{2}+\beta(y).
\end{array}\right.
\end{equation}

\noindent Finally, we recall the following result from \cite{sck2010}, which states the existence and uniqueness of regular enough kernel functions involved in the backstepping transformation.

\begin{thm}[Theorem 3.1 in \cite{sck2010}]
\label{kernel}
Let $\lambda\in C^{2}([0,L])$ and $\beta\in C([0,L])$. Then, there exist a unique $k=k(x,y)$ in $C^{2}(\mathcal{T})$ and a unique $s=s(x,y)$ in $C^{2}(\mathcal{T})$ such that

\begin{equation}
\label{kernel_1}
\left\{\begin{array}{rl}
k_{xx}-k_{yy}=\rho_{1}s_{yy}+\rho_{2}k+\rho_{3}s+\rho_{4}s_{y},&~(x,y)\in\mathcal{T},\\[1ex]
k(x,x)=m(x),&~x\in[0,L],\\[1ex]
k(x,0)=0,&~x\in[0,L],
\end{array}\right.
\end{equation}

\begin{equation}
\label{kernel_2}
\left\{\begin{array}{rl}
s_{xx}-s_{yy}=\rho_{1}k+\rho_{5}s,&~(x,y)\in\mathcal{T},\\[1ex]
s(x,x)=-\mbox{sinh}\left(\int_{0}^{x}a(\tau)\,d\tau\right),&~x\in[0,L],\\[1ex]
s(x,0)=0,&~x\in[0,L].
\end{array}\right.
\end{equation}
\end{thm}

We have introduced all the required elements. We continue with the feedback design. In \eqref{dd_1} we consider $U(t)=U_{1}(t)+U_{2}(t)$. Then, in virtue of the computations done in \cite[Section 3.1]{sck2010}, we employ the backstepping transformation

\begin{equation}
\label{BS}
w(t,x)=h(x)u(t,x)-\int_{0}^{x}k(x,y)u(t,y)\,dy-\int_{0}^{x}s(x,y)u_{t}(t,y)\,dy
\end{equation}

\noindent and choose

\begin{equation}
\label{control_1}
U_{1}(t)=\frac{1}{h(L)}\left(\int_{0}^{L}k(L,y)u(t,y)\,dy+\int_{0}^{L}s(L,y)u_{t}(t,y)\,dy\right)
\end{equation}

\noindent to transform \eqref{dd_1} into the target system

\begin{equation}
\label{dd_2}
\left\{\begin{array}{rl}
w_{tt}=w_{xx}-2dw_{t}-d^{2}w,&~(t,x)\in(0,\infty)\times(0,L),\\[1ex]
w(t,0)=0,~w(t,L)=h(L)p(t)+h(L)U_{2}(t),&~t\in(0,\infty),\\[1ex]
w(0,x)=w_{0}(x),~w_{t}(0,x)=w_{1}(x),&~x\in(0,L),
\end{array}\right.
\end{equation}

\noindent where the transformed initial condition $(w_{0},w_{1})$, as stated in \cite[Section 3.3]{sck2010}, is given by

\begin{equation}
\label{transformed_initial}
\arraycolsep=2pt
\left\{\begin{array}{lll}
w_{0}(x) & = & h(x)u_{0}(x)-\int_{0}^{x}k(x,y)u_{0}(y)\,dy-\int_{0}^{x}s(x,y)u_{1}(y)\,dy,\\[2ex]
w_{1}(x) & = & h(x)u_{1}(x)+s_{y}(x,x)u_{0}(x)-s(x,x)u_{0}'(x)\\[1ex]
 & & -\int_{0}^{x}\left[2\lambda(y)s(x,y)+k(x,y)\right]u_{1}(y)\,dy-\int_{0}^{x}\left[\beta(y)s(x,y)+s_{yy}(x,y)\right]u_{0}(y)\,dy.
\end{array}\right.
\end{equation}

As an intermediate step, we introduce a modified energy as follows: for a regular enough state $(w,w_{t})$ of the target system \eqref{dd_2}, we define its modified energy by

\begin{equation}
\label{m_energy}
V(t)=\frac{1}{2}\int_{0}^{L}(w_{t}+dw)^{2}\,dx+\frac{1}{2}\int_{0}^{L}|w_{x}|^{2}\,dx,~t\in[0,\infty).
\end{equation}

We proceed to perform formal computations, which later are justified in view of the analysis done in Section \ref{dd_wp}. Performing integration by parts when needed, we have

\begin{equation*}
\begin{split}
V'(t) =& \int_{0}^{L}(w_{t}+dw)(w_{tt}+dw_{t})\,dx+\int_{0}^{L}w_{x}w_{xt}\,dx\\
 =& \int_{0}^{L}\left(w_{t}(w_{xx}-2dw_{t}-d^{2}w)+d|w_{t}|^{2}+dw(w_{xx}-2dw_{t}-d^{2}w)+d^{2}ww_{t}\right)dx\\
 &-\int_{0}^{L}w_{xx}w_{t}\,dx++w_{x}(t,L)w_{t}(t,L)\\
 =& \int_{0}^{L}\left(-d|w_{t}|^{2}+dww_{xx}-2d^{2}ww_{t}-d^{3}|w|^{2}\right)dx+w_{x}(t,L)w_{t}(t,L)\\
 =&-2d\int_{0}^{L}\frac{1}{2}\left(|w_{t}|^{2}+2dww_{t}+d^{2}|w|^{2}+|w_{x}|^{2}\right)dx+dw(t,L)w_{x}(t,L)+w_{x}(t,L)w_{t}(t,L).
\end{split}
\end{equation*}

\noindent Accordingly, we arrive at

\begin{equation}
\label{m_energy_1}
V'(t)=-2dV(t)+w_{x}(t,L)\left(w_{t}(t,L)+dw(t,L)\right).
\end{equation}

From the previous expression we shall design $U_{2}(t)$ as a feedback law. Let us assume \textbf{(A1)}. Then, as $w(t,L)=h(L)p(t)+h(L)U_{2}(t)$ it follows

\begin{equation}
\label{m_energy_2}
V'(t)+2dV(t)=w_{x}(t,L)w_{t}(t,L)+w_{x}(t,L)dh(L)p(t)+w_{x}(t,L)dh(L)U_{2}(t).
\end{equation}

\noindent Then, with the choice of $U_{2}(t)$ so that

\begin{equation}
\label{m_energy_3}
U_{2}(t)=-\frac{1}{dh(L)}w_{t}(t,L)-P\,\text{sign}\left(w_{x}(t,L)\right),
\end{equation}

\noindent the right-hand side of \eqref{m_energy_2} becomes non-positive. Indeed, thanks to the property $\theta f=|f|$ for every $\theta\in\text{sign}(f)$, we have

\begin{align}
\notag
 &w_{x}(t,L)w_{t}(t,L)+w_{x}(t,L)dh(L)p(t)+w_{x}(t,L)dh(L)U_{2}(t)\\[1ex]
\notag
 &\leq |w_{x}(t,L)|dh(L)P+w_{x}(t,L)dh(L)\theta=|w_{x}(t,L)|dh(L)P-|w_{x}(t,L)|dh(L)P\\[1ex]
\label{m_energy_4}
 &=0\mbox{ for every }\theta\in-P\,\text{sign}\left(w_{x}(t,L)\right).
\end{align}

\newpage

\noindent Therefore, in virtue of \eqref{m_energy_2} and \eqref{m_energy_4} we conclude that the solutions of the corresponding closed-loop system, obtained by considering the target system \eqref{dd_2} together with \eqref{m_energy_3}, satisfy

\begin{equation}
\label{dd_3}
V(t)\leq e^{-2dt}V(0),~t\in[0,\infty).
\end{equation}

The closed-loop system is the following one.

\begin{equation}
\label{dd_4}
\left\{\begin{array}{rl}
w_{tt}=w_{xx}-2dw_{t}-d^{2}w,&~(t,x)\in(0,\infty)\times(0,L),\\[1ex]
w(t,0)=0,&~t\in(0,\infty),\\[1ex]
w(t,L)+\frac{1}{d}w_{t}(t,L)\in h(L)p(t)-h(L)P\,\text{sign}\left(w_{x}(t,L)\right),&~t\in(0,\infty),\\[1ex]
w(0,x)=w_{0}(x),~w_{t}(0,x)=w_{1}(x),&~x\in(0,L).
\end{array}\right.
\end{equation}


The next result states an equivalence between the modified energy and the energy.

\begin{lem}
\label{e_energy}
Let $(w,w_{t})$ be such that $w(t,\cdot)\in H^{1}(0,L)$, $w(t,0)=0$ and $w_{t}(t,\cdot)\in L^{2}(0,L)$ for every $t\in[0,\infty)$. Then,

\begin{align}
\notag
 & \mbox{min}\left\{\frac{1}{8L^{2}d^{2}+1},\frac{1}{2}\right\}\int_{0}^{L}(|w_{t}|^{2}+|w_{x}|^{2})\,dx\leq 2V(t)\\
\label{e_energy_1}
 & \leq\mbox{max}\left\{8L^{2}d^{2}+1,2\right\}\int_{0}^{L}(|w_{t}|^{2}+|w_{x}|^{2})\,dx,~t\in[0,\infty).
\end{align}
\end{lem}

\noindent\textbf{Proof~~} By hypotheses it is available Poincar\'e inequality $\int_{0}^{L}|w(t,x)|^{2}\,dx\leq 4L^{2}\int_{0}^{L}|w_{x}(t,x)|^{2}\,dx$ for every $t\in[0,\infty)$. We prove the first inequality only (from left to right), since the proof of the second one is similar. Thanks to Poincar\'e and Cauchy inequalities, for every $\delta>0$ we have

\begin{align}
\notag
2V(t) &= \int_{0}^{L}\left(|w_{t}|^{2}+2w_{t}dw+d^{2}|w|^{2}\right)\,dx+\int_{0}^{L}|w_{x}|^{2}\left(\frac{1}{2}+\frac{1}{2}\right)\,dx\\
\label{e_energy_2}
 &\geq \int_{0}^{L}|w_{t}|^{2}(1-d\delta)\,dx+\int_{0}^{L}|w_{x}|^{2}\frac{1}{2}\,dx+\int_{0}^{L}|w|^{2}\left(d^{2}-\frac{d}{\delta}+\frac{1}{8L^{2}}\right)\,dx.
\end{align}

\noindent By chosing $\delta>0$ so that $d^{2}-\frac{d}{\delta}+\frac{1}{8L^{2}}=0$, which is when $\delta=\frac{8L^{2}d}{8L^{2}d^{2}+1}$, it follows that $1-d\delta=\frac{1}{8L^{2}d^{2}+1}$. Accordingly, from \eqref{e_energy_2} we get the first inequality because

\begin{equation*}
2V(t)\geq\frac{1}{8L^{2}d^{2}+1}\int_{0}^{L}|w_{t}|^{2}\,dx+\frac{1}{2}\int_{0}^{L}|w_{x}|^{2}\,dx.
\end{equation*}

The second inequality also follows by applying Poincar\'e and Cauchy inequalities.\hfill$\blacksquare$\\

To conclude the feedback design part of Theorem \ref{dd_mr}, let us obtain \eqref{rapid} for \eqref{dd_1}. To that end, we use that both the backstepping operator and its inverse are linear and bounded operators. Let us introduce the space

\begin{equation}
\label{spaces}
H=\{h\in H^{1}(0,L)~/~h(0)=0\}\mbox{ and }X=H\times L^{2}(0,L),
\end{equation}

\noindent with norm $\|(a,b)\|_{X}^{2}=\|b\|_{L^{2}(0,L)}^{2}+\|a'\|_{L^{2}(0,L)}^{2}$. From \cite[Section 3.3]{sck2010} we have that the backstepping operator $\Pi:X\to X$ given by $\Pi(u_{0},u_{1})=(w_{0},w_{1})$, with $(w_{0},w_{1})$ defined in \eqref{transformed_initial}, has a linear and bounded inverse $\Pi^{-1}:X\to X$.

\newpage
\mbox{}
\vspace{-1.25cm}

\begin{rem}
\label{BS_operator}
As stated in \cite[Section 3.3]{sck2010}, the backstepping operator $\Pi$ maps the state $(u,u_{t})$ of \eqref{dd_1} with \eqref{control_1} into the state $(w,w_{t})$ of the target system \eqref{dd_2}.
\end{rem}

Let us set $C_{1}=\|\Pi^{-1}\|_{\mathcal{L}(X)}$ and $C_{2}=\|\Pi\|_{\mathcal{L}(X)}$, which satisfy $C_{1}C_{2}\geq1$. As the energy satisfies $2E(t)=\|(u,u_{t})(t,\cdot)\|_{X}^{2}$ in virtue of \eqref{energy}, it follows from the continuity of $\Pi^{-1}$, \eqref{e_energy_1}, \eqref{dd_3} and the continuity of $\Pi$ that

\begin{align}
\notag
\|(u,u_{t})(t,\cdot)\|_{X}^{2} & \leq C_{1}^{2}\|(w,w_{t})(t,\cdot)\|_{X}^{2}\\[2ex]
\notag
 & \leq C_{1}^{2}\frac{1}{\mbox{min}\left\{\frac{1}{8L^{2}d^{2}+1},\frac{1}{2}\right\}}2V(t)\\[2ex]
\notag
 & \leq C_{1}^{2}\frac{1}{\mbox{min}\left\{\frac{1}{8L^{2}d^{2}+1},\frac{1}{2}\right\}}e^{-2dt}2V(0)\\[2ex]
\notag
 & \leq C_{1}^{2}\frac{\mbox{max}\left\{8L^{2}d^{2}+1,2\right\}}{\mbox{min}\left\{\frac{1}{8L^{2}d^{2}+1},\frac{1}{2}\right\}}e^{-2dt}\|(w_{0},w_{1})\|_{X}^{2}\\[2ex]
\label{decay}
 & \leq C_{1}^{2}C_{2}^{2}\frac{\mbox{max}\left\{8L^{2}d^{2}+1,2\right\}}{\mbox{min}\left\{\frac{1}{8L^{2}d^{2}+1},\frac{1}{2}\right\}}e^{-2dt}\|(u_{0},u_{1})\|_{X}^{2}\mbox{ for all }t\in[0,\infty).
\end{align}

Accordingly, we have shown \eqref{rapid} for \eqref{dd_1}. In other words, we have shown the feedback design part of Theorem \ref{dd_mr} and \eqref{dd_mr_decay} with $C=C_{1}C_{2}\sqrt{K}$, where

\begin{equation}
\label{constant}
K=\mbox{max}\left\{8L^{2}d^{2}+1,2\right\}/\mbox{min}\left\{\frac{1}{8L^{2}d^{2}+1},\frac{1}{2}\right\}.
\end{equation}

\begin{rem}
As stated in \cite[Section 2.6]{vabk2026}, let us note that the inverse backstepping transformation is not used in the feedback design. Nevertheless, its existence and properties are indeed used for the proof of the closed-loop stability of \eqref{dd_mr_closed}, as we have just exposed for the obtention of \eqref{decay}.
\end{rem}

\subsection{Well-posedness}
\label{dd_wp}

\hspace{0.75cm}In this section we apply the maximal monotone operator theory to prove the well-posedness of the closed-loop system \eqref{dd_4}, which is a differential inclusion, with the inclusion being located at the\l boundary condition.\\

Let us recall \eqref{spaces}. Also, on the Hilbert space $X$ let us consider the inner product

\begin{equation}\label{inn-prod}
    \langle (q_1,l_1), (q_2,l_2) \rangle_{X} =
    \int_0^L q_1'(x)q_2'(x)\,dx + \int_0^L  (l_1(x)+dq_1(x))(l_2(x)+dq_2(x))\,dx.
\end{equation}

\noindent Let us introduce

\begin{equation}
\label{lift_dd}
y(t,x)=w(t,x)-\psi(t)\phi(x)~\mbox{with}~\psi(t)=dh(L)\left(e^{-dt}\int_{0}^{t}e^{ds}p(s)\,ds\right)~\mbox{and}~\phi(x)=\frac{x}{L}\left(2-\frac{x}{L}\right).
\end{equation}

\noindent Since $\phi(0)=0$, $\phi(L)=1$ and $\phi'(L)=0$ we get that $y(t,0)=w(t,0)$, $y(t,L)+\frac{1}{d}y_{t}(t,L)=w(t,L)+\frac{1}{d}w_{t}(t,L)-h(L)p(t)$ and $y_{x}(t,L)=w_{x}(t,L)$. Let us assume \textbf{(A2)$\mathbf{_{DD}}$}, thus obtaining that $y(0,x)=w(0,x)$ and $y_{t}(0,x)=w_{t}(0,x)$. Then, it follows from \eqref{dd_4} that $y=y(t,x)$ satisfies

\begin{equation}
\label{dd_5}
\left\{\begin{array}{rl}
y_{tt}=y_{xx}-2dy_{t}-d^{2}y+f,&~(t,x)\in(0,\infty)\times(0,L),\\[1ex]
y(t,0)=0,&~t\in(0,\infty),\\[1ex]
y(t,L)+\frac{1}{d}y_{t}(t,L)\in -h(L)P\,\text{sign}\left(y_{x}(t,L)\right),&~t\in(0,\infty),\\[1ex]
y(0,x)=w_{0}(x),~y_{t}(0,x)=w_{1}(x),&~x\in(0,L),
\end{array}\right.
\end{equation}

\noindent where $f(x,t) = (-\partial_{tt} +\partial_{xx} - 2d \partial_t - d^2)(\psi(t)\phi(x))$. Finally, let us consider the operator

\begin{equation}
\label{operator_dd}
\left\{\begin{array}{l}
\mathcal{A}:D(\mathcal{A})\subset X\to X,\\[2ex]
\mathcal{A}(q,l)=(-l,-q''+2dl+d^{2}q),\\[2ex]
D(\mathcal{A})=\{(q,l)\in X~/~\mathcal{A}(q,l)\in X,~dq(L)+l(L)\in-dh(L)P\,\text{sign}(q'(L))\}.
\end{array}\right.
\end{equation}

\begin{rem}
\label{ope_rem_dd}
We have that $D(\mathcal{A})$ is a non-empty set because $H_{0}^{2}(0,L)\times H_{0}^{1}(0,L)\subset D(\mathcal{A})$ due to the true statement $0\in-dh(L)P[-1,1]$. Because of the sign multivalued operator \eqref{sign} it follows that $D(\mathcal{A})$ is not a linear subspace, and hence, $\mathcal{A}$ is not a linear operator. Also, $D(\mathcal{A})$ is dense in $X$. Finally,\\

$D(\mathcal{A})=\{(q,l)\in(H^{2}(0,L)\cap H)\times H~/~dq(L)+l(L)\in-dh(L)P\,\text{sign}(q'(L))\}$.
\end{rem}

Thanks to \eqref{operator_dd} we can write \eqref{dd_5} in operator form as follows:

\begin{equation}
\label{operator_dd_2}
\left\{\begin{array}{rl}
\dfrac{d}{dt}(y,y_{t})+\mathcal{A}(y,y_{t})=(0,f),&t\in[0,\infty),\\[2ex]
(y,y_{t})(0)=(w_{0},w_{1}).
\end{array}\right.
\end{equation}

We proceed to prove that \eqref{operator_dd_2} is well-posed by applying the maximal monotone operator theory. In that direction we have the following two results: the first one (Proposition \ref{wp1_dd}) states that the operator $\mathcal{A}$ is monotone, while the second one (Proposition \ref{wp2_dd}) states that the operator $I+\mathcal{A}$ is surjective.

\begin{prop}
\label{wp1_dd}
The operator defined by \eqref{operator_dd} is monotone.
\end{prop}

\noindent\textbf{Proof~~}Let $(q_1,l_1)$ and $(q_2,l_2)$ be in $\mathcal{D}(\mathcal{A})$. For brevity in notation and only in this proof, let us set $q= q_1-q_2$ and $l=l_1-l_2$. Then, in virtue of the inner product \eqref{inn-prod}, we have, after two integrations by parts,

\begin{align*}
 &\langle \mathcal{A} (q_1,l_1) - \mathcal{A}(q_2,l_2), (q_1,l_1)-(q_2,l_2) \rangle_{X}\\[1ex]
 &=\int_{0}^{L}(-l)'q'\,dx+\int_{0}^{L}\left((-q''+2dl+d^{2}q)+d(-l)\right)(l+dq)\,dx\\[1ex]
 &=d\|q'\|_{L^{2}(0,L)}^{2}+d^{3}\|q\|_{L^{2}(0,L)}^{2}+d\|l\|_{L^{2}(0,L)}^{2}+2d^{2}\int_{0}^{L}ql\,dx+\left(-l(L)q'(L)-dq(L)q'(L)\right)\\[1ex]
 &= I_{1}+I_{2}.
\end{align*}

On the one hand, since $\left|2d^{2}\int_{0}^{L}ql\,dx\right|\leq d^{3}\|q\|_{L^{2}(0,L)}^{2}+d\|l\|_{L^{2}(0,L)}^{2}$ we get that $I_{1}\geq d\|q'\|_{L^{2}(0,L)}^{2}\geq0$. On the other hand, as $dq_{j}(L)+l_{j}(L)\in-dh(L)P\,\text{sign}(q_{j}'(L))$ for $j\in\{1,2\}$, there exists $\theta_{j}\in\text{sign}(q_{j}'(L))$ such that $dq_{j}(L)+l_{j}(L)=-dh(L)P\theta_{j}$. Therefore,

\begin{align*}
 &I_{2}=-l(L)q'(L)-dq(L)q'(L)\\[1ex]
 &=-\left((l_{1}(L)+dq_{1}(L))-(l_{2}(L)+dq_{2}(L))\right)(q_{1}'(L)-q_{2}'(L))\\[1ex]
 &=dh(L)P(\theta_{1}-\theta_{2})(q_{1}'(L)-q_{2}'(L)).
\end{align*}

\noindent We conclude that $I_{2}\geq0$ because the sign multivalued operator \eqref{sign} is monotone. Accordingly, the operator $\mathcal{A}$ is monotone as $I_{1}+I_{2}\geq0$.\hfill$\blacksquare$

\begin{prop}
\label{wp2_dd}
The operator defined by \eqref{operator_dd} satisfies $R(I+\mathcal{A})=X$.
\end{prop}

\noindent\textbf{Proof~~}Given $(m,n)\in X$ we need to prove the existence of $(q,l)\in D(\mathcal{A})$, domain being defined in \eqref{operator_dd}, such that $(q,l)+\mathcal{A}(q,l)=(m,n)$. To that end, let us introduce the problem

\begin{equation}
\label{wp2_dd_1}
\left\{\begin{array}{l}
-q''+q(1+2d+d^{2})=n+m(1+2d),\\[1ex]
q(0)=0,~dq(L)+(q(L)-m(L))\in-dh(L)P\,\text{sign}(q'(L)).
\end{array}\right.
\end{equation}

\noindent By proving that \eqref{wp2_dd_1} possesses a unique solution $q\in H^{2}(0,L)$ and by putting $l=q-m$, then we get that $(q,l)\in D(\mathcal{A})$ and that $(q,l)+\mathcal{A}(q,l)=(m,n)$. For the proof we proceed as in \cite[Section 2]{cp1990} and analyze an optimization problem.

\begin{rem}
The analysis in \cite[Section 2]{cp1990} has also been adapted for the proof of \cite[Proposition 3.1]{gh2023} and the proof of \cite[Proposition 3.2]{gh2025}.
\end{rem}

Let us introduce the functional $J:H^{2}(0,L)\cap H\to\mathbb{R}$ by

\begin{align}
\notag
J(q)= &~\frac{1}{2}\int_{0}^{L}\left[(q'')^{2}+(q')^{2}(1+2d+d^{2})+2q''\left(n+m(1+2d)\right)\right]dx\\[1ex]
\label{functional_dd}
&+(1+d)\left[-m(L)q'(L)+dh(L)P\,\varphi_{\sigma}(q'(L))\right],
\end{align}

\noindent where $\varphi_{\sigma}:\mathbb{R}\to\mathbb{R}$ is the Moreau regularization of the convex and continuous function $\varphi:\mathbb{R}\to\mathbb{R}$ given by $\varphi(x)=|x|$. Thanks to \cite[Chapter IV, Proposition 1.8]{s1997} (see also \cite[Theorem 2.9]{b2010}) we have

\begin{equation*}
\varphi_{\sigma}(x)=\frac{\sigma}{2}|\alpha_{\sigma}(x)|^{2}+\varphi(J_{\sigma}(x)),~(\sigma,x)\in(0,\infty)\times\mathbb{R},
\end{equation*}

\noindent where $\alpha_{\sigma}=\sigma^{-1}(I-J_{\sigma}):\mathbb{R}\to\mathbb{R}$ is the Yosida approximation of the maximal monotone operator $\alpha:\mathbb{R}\to 2^{\mathbb{R}}$ given by $\alpha(x)=(\partial\varphi)(x)=\text{sign}(x)$ and $J_{\sigma}=(I+\sigma\alpha)^{-1}:\mathbb{R}\to\mathbb{R}$ is the resolvent of $\alpha$. Furthermore, also from \cite[Chapter IV, Proposition 1.8]{s1997} (see also \cite[Theorem 2.9]{b2010}) we have that $\varphi_{\sigma}$ is a convex and differentiable function satisfying that $\varphi_{\sigma}'(x)=\alpha_{\sigma}(x)$.\\

Since $q\in H^{2}(0,L)\cap H$ we have the following inequalities.

\begin{equation}
\label{ineq_wp}
\|q\|_{L^{2}(0,L)}^{2}\leq 4L^{2}\|q'\|_{L^{2}(0,L)}^{2},~~ |q(L)|\leq2L^{1/2}\|q'\|_{L^{2}(0,L)}~~\mbox{and}~~|q'(L)|\leq(1+L^{-1})^{1/2}\|q\|_{H^{2}(0,L)}.
\end{equation}

The functional $J$ is convex and continuous, the latter being a consequence of the continuous injection of $H^{2}(0,L)$ into $C^{1}([0,L])$. Also, considering the non-negativity of $\varphi_{\sigma}$ in \eqref{functional_dd} and applying \eqref{ineq_wp}, there exist two positive constants, $C_{1}=C_{1}(L)$ and $C_{2}=C_{2}(L,d,m,n)$, such that

\begin{equation*}
J(q)\geq C_{1}\|q\|_{H^{2}(0,L)}^{2}-C_{2}\|q\|_{H^{2}(0,L)},
\end{equation*}

\noindent from which it follows that $J(q)\to\infty$ as $\|q\|_{H^{2}(0,L)}\to\infty$, concluding that the functional $J$ is coercive. Therefore, \cite[Theorem 2.19]{p2015} yields the existence of a minimizer $q\in H^{2}(0,L)\cap H$ for the functional $J$. Then, in virtue of \cite[Proposition 3.20]{p2015} and \cite[Theorem 3.24]{p2015}, we have that the G\^ateaux derivative of the functional $J$ at such minimizer must vanish for each direction $\rho\in H^{2}(0,L)\cap H$. Accordingly, as $\varphi_{\sigma}'(x)=\alpha_{\sigma}(x)$ we obtain

\begin{align*}
J'(q;\rho)= &~\int_{0}^{L}\left[q''\rho''+q'\rho'(1+2d+d^{2})+\rho''\left(n+m(1+2d)\right)\right]dx\\[1ex]
&+(1+d)\left[-m(L)+dh(L)P\,\alpha_{\sigma}(q'(L))\right]\rho'(L)=0~~\forall \rho\in H^{2}(0,L)\cap H,
\end{align*}

\noindent and after one integration,

\begin{align}
\notag
J'(q;\rho)= &~\int_{0}^{L}\left[q''-q(1+2d+d^{2})+\left(n+m(1+2d)\right)\right]\rho''\,dx\\[1ex]
\notag
&+(1+d)\left[dq(L)+(q(L)-m(L))+dh(L)P\,\alpha_{\sigma}(q'(L))\right]\rho'(L)=0\\[2ex]
\label{derivative_dd}
&\forall\rho\in H^{2}(0,L)\cap H.
\end{align}

By considering $C_{0}^{\infty}(0,L)\subset H^{2}(0,L)\cap H$ in \eqref{derivative_dd}, we can deduce thanks to \cite[Lemma 8.1]{bre2010} the existence of two constants, $C_{3}\in\mathbb{R}$ and $C_{4}\in\mathbb{R}$, so that $q''-q(1+2d+d^{2})+\left(n+m(1+2d)\right)=C_{3}x+C_{4}$ for almost every $x\in(0,L)$. Then, from \eqref{derivative_dd} we get

\begin{align}
\notag
&-C_{3}\rho(L)+C_{3}L\rho'(L)+C_{4}(\rho'(L)-\rho'(0))\\[1ex]
\label{wp2_dd_2}
&+(1+d)\left[dq(L)+(q(L)-m(L))+dh(L)P\,\alpha_{\sigma}(q'(L))\right]\rho'(L)=0~~\forall\rho\in H^{2}(0,L)\cap H.
\end{align}

\noindent From \eqref{wp2_dd_2} we can deduce that $C_{3}=0$ by chosing a third degree polynomial satisfying $\rho(L)=1$ and $\rho(0)=\rho'(0)=\rho'(L)=0$. Similarly, we can deduce that $C_{4}=0$ by chosing a second degree polynomial satisfying $\rho'(0)=1$ and $\rho(0)=\rho'(L)=0$.

\begin{rem}
Those two polynomials do exist. Indeed, the conclusion $C_{3}=0$ can be obtained with\l $\rho(x)=-2L^{-3}x^{3}+3L^{-2}x^{2}$, and $C_{4}=0$, with $\rho(x)=-(2L)^{-1}x^{2}+x$. Both polynomials satisfy the requirements.
\end{rem}

\noindent Finally, we can also deduce that

\begin{equation*}
dq(L)+(q(L)-m(L))+dh(L)P\,\alpha_{\sigma}(q'(L))=0.
\end{equation*}

Accordingly, so far we have shown:

\begin{lem}
\label{wp3_dd}
For any $\sigma\in(0,\infty)$ there exists $q_{\sigma}\in H^{2}(0,L)$ such that $-q_{\sigma}''+q_{\sigma}(1+2d+d^{2})=n+m(1+2d)$ for almost every $x\in(0,L)$, $q_{\sigma}(0)=0$ and $dq_{\sigma}(L)+(q_{\sigma}(L)-m(L))=-dh(L)P\,\alpha_{\sigma}(q_{\sigma}'(L))$.
\end{lem}

In view of Lemma \ref{wp3_dd} it seems that we are close to prove that \eqref{wp2_dd_1} possesses a unique solution $q\in H^{2}(0,L)$. To that end, let $q_{\sigma}$ be given by Lemma \ref{wp3_dd} and let us see what happens to it when $\sigma\to0^{+}$. We proceed to obtain upper bounds, independent of $\sigma$, for $\|q_{\sigma}\|_{H^{2}(0,L)}$ and $|\alpha_{\sigma}(q_{\sigma}'(L))|$.

\newpage

$\bullet$ \underline{Upper bound for $\|q_{\sigma}\|_{H^{2}(0,L)}$ independent of $\sigma$}\\

In virtue of Lemma \ref{wp3_dd} we can get

\begin{equation*}
\|q_{\sigma}'\|_{L^{2}(0,L)}^{2}-q_{\sigma}'(L)q_{\sigma}(L)+\frac{1}{2}\|q_{\sigma}\|_{L^{2}(0,L)}^{2}\leq\frac{1}{2}\|n+m(1+2d)\|_{L^{2}(0,L)}^{2}.
\end{equation*}

\noindent In order to treat the term $-q_{\sigma}'(L)q_{\sigma}(L)$, let us note that from $0\in\alpha(0)$ we conclude that $J_{\sigma}(0)=0$, and consequently, $\alpha_{\sigma}(0)=0$. Then, by Lemma \ref{wp3_dd} and as $\alpha_{\sigma}$ is a monotone operator, we infer

\begin{equation*}
-q_{\sigma}'(L)q_{\sigma}(L)=-\frac{1}{1+d}m(L)q_{\sigma}'(L)+\frac{d}{1+d}h(L)P\alpha_{\sigma}(q_{\sigma}'(L))q_{\sigma}'(L)\geq-\frac{1}{1+d}m(L)q_{\sigma}'(L).
\end{equation*}

\noindent Therefore, applying \eqref{ineq_wp} and Cauchy inequality, for every $\delta>0$ we have

\begin{equation}
\label{uniform_dd_1}
\|q_{\sigma}'\|_{L^{2}(0,L)}^{2}+\frac{1}{2}\|q_{\sigma}\|_{L^{2}(0,L)}^{2}\leq\frac{1}{2}\|n+m(1+2d)\|_{L^{2}(0,L)}^{2}+(1+L^{-1})\frac{|m(L)|^{2}}{(1+d)^{2}}\frac{\delta}{2}+\frac{1}{2\delta}\|q_{\sigma}\|_{H^{2}(0,L)}^{2}.
\end{equation}

Thanks to Lemma \ref{wp3_dd} we can get

\begin{equation}
\label{uniform_dd_2}
\|q_{\sigma}''\|_{L^{2}(0,L)}^{2}\leq2(1+d)^{4}\|q_{\sigma}\|_{L^{2}(0,L)}^{2}+2\|n+m(1+2d)\|_{L^{2}(0,L)}^{2}.
\end{equation}

Finally, from the combination of \eqref{uniform_dd_1} together with \eqref{uniform_dd_2} and a suitable choice for $\delta>0$, we can obtain the existence of two positive constants, $C_{5}=C_{5}(d)$ and $C_{6}=C_{6}(L,d)$, such that

\begin{equation}
\label{uniform_dd_3}
\|q_{\sigma}\|_{H^{2}(0,L)}^{2}\leq C_{5}\|n+m(1+2d)\|_{L^{2}(0,L)}^{2}+C_{6}|m(L)|^{2}.
\end{equation}

$\bullet$ \underline{Upper bound for $|\alpha_{\sigma}(q_{\sigma}'(L))|$ independent of $\sigma$}\\

By Lemma \ref{wp3_dd} it follows that

\begin{equation*}
|\alpha_{\sigma}(q_{\sigma}'(L))|\leq\frac{(d+1)}{dh(L)P}|q(L)|+\frac{1}{dh(L)P}|m(L)|.
\end{equation*}

\noindent Then, in view of \eqref{ineq_wp} and \eqref{uniform_dd_3}, we see that there exist two positive constants, $C_{7}=C_{7}(L,d)$ and $C_{8}=C_{8}(L,d)$, such that

\begin{equation}
\label{uniform_dd_4}
|\alpha_{\sigma}(q_{\sigma}'(L))|\leq\frac{C_{7}}{h(L)P}\|n+m(1+2d)\|_{L^{2}(0,L)}+\frac{C_{8}}{h(L)P}|m(L)|.
\end{equation}

Therefore, in virtue of \eqref{uniform_dd_3} and \eqref{uniform_dd_4} we conclude that the sequences $\left\{q_{\sigma}\right\}_{\sigma>0}\subset H^{2}(0,L)$ and $\left\{\alpha_{\sigma}(q_{\sigma}'(L))\right\}_{\sigma>0}\subset\mathbb{R}$ are bounded. Then, there exist $(q,b)\in H^{2}(0,L)\times\mathbb{R}$ and subsequences, which we denote by the same symbols, such that $q_{\sigma}\rightharpoonup q$ in $H^{2}(0,L)$ (weak convergence) and $\alpha_{\sigma}(q_{\sigma}'(L))\to b$ in $\mathbb{R}$ as $\sigma\to0^{+}$. Moreover, we infer that $q_{\sigma}\to q$ in $C^{1}([0,L])$ as $\sigma\to0^{+}$ since the injection of $H^{2}(0,L)$ into $C^{1}([0,L])$ is compact, implying that $q_{\sigma}(x)\to q(x)$ and $q_{\sigma}'(x)\to q'(x)$ in $\mathbb{R}$ for all $x\in[0,L]$ as $\sigma\to0^{+}$.\\

Accordingly, Lemma \ref{wp3_dd} and the previous arguments yield:

\begin{lem}
\label{wp4_dd}
There exists $q\in H^{2}(0,L)$ such that $-q''+q(1+2d+d^{2})=n+m(1+2d)$ for almost every $x\in(0,L)$, $q(0)=0$ and $dq(L)+(q(L)-m(L))=-dh(L)P\,b$.
\end{lem}

\newpage

In view of Lemma \ref{wp4_dd} and \eqref{wp2_dd_1} we see that in order to complete the proof of $R(I+\mathcal{A})=X$ we just need to prove that $b\in\alpha(q'(L))=\text{sign}(q'(L))$. The required arguments can be found at the end of the proof of \cite[Proposition 3.1]{gh2023} or at the end of the proof of \cite[Proposition 3.2]{gh2025}. The proof of Proposition \ref{wp2_dd} is complete.\hfill$\blacksquare$\\

Let us return to \eqref{operator_dd_2}. Due to Proposition \ref{wp1_dd} and Proposition \ref{wp2_dd} we can apply \cite[Chapter IV, Lemma 1.3]{s1997} to conclude that the operator $\mathcal{A}$ defined in \eqref{operator_dd} is maximal monotone. Regarding the data $(0,f)$ and $(w_{0},w_{1})$ in \eqref{operator_dd_2}: recalling that $f(x,t) = (-\partial_{tt} +\partial_{xx} - 2d \partial_t - d^2)(\psi(t)\phi(x))$, by \textbf{(A2)$\mathbf{_{DD}}$} we get that $f\in W^{1,1}(0,\infty;C^{\infty}(\mathbb{R}))$; also, let us assume that $(w_{0},w_{1})\in D(\mathcal{A})$. Then, we can apply\l \cite[Chapter IV, Theorem 4.1]{s1997} to conclude the existence of a unique $(y,y_{t})\in W^{1,1}(0,\infty;X)$ such that

\begin{equation}
\label{wp5_dd}
\left\{\begin{array}{l}
\dfrac{d}{dt}(y,y_{t})+\mathcal{A}(y,y_{t})=(0,f)\mbox{ for almost every }t>0,\\[2ex]
(y,y_{t})\in D(\mathcal{A})\mbox{ for every }t\geq0,\\[2ex]
(y,y_{t})(0)=(w_{0},w_{1}).
\end{array}\right.
\end{equation}

With \eqref{dd_4} (the closed-loop system), \eqref{lift_dd} and \eqref{dd_5} in mind, due to \eqref{wp5_dd} we conclude that $(w,w_{t})\in W^{1,1}(0,\infty;X)$ uniquely solves

\begin{equation}
\label{wp6_dd}
\arraycolsep=0pt
\left\{\begin{array}{l}
w_{tt}=w_{xx}-2dw_{t}-d^{2}w ~\mbox{for almost every}~(t,x)\in(0,\infty)\times(0,L),\\[1ex]
w(t,0)=0~\mbox{for every}~t\in[0,\infty),\\[1ex]
w(t,L)+\frac{1}{d}w_{t}(t,L)\in h(L)p(t)-h(L)P\,\text{sign}\left(w_{x}(t,L)\right)~\mbox{for every}~t\in[0,\infty),\\[1ex]
w(0,x)=w_{0}(x)~\mbox{and}~w_{t}(0,x)=w_{1}(x)~\mbox{for every}~x\in(0,L).
\end{array}\right.
\end{equation}

We have shown the well-posedness part of Theorem \ref{dd_mr}.

\subsection{Conclusion}
\label{dd_c}

\hspace{0.75cm}Here we gather the conclusions of both the feedback design (see Section \ref{dd_fd}) and the well-posedness (see Section \ref{dd_wp}) parts to establish our main result regarding the Dirichlet-Dirichlet case.\\

In Section \ref{dd_fd} we mentioned that from \cite[Section 3.3]{sck2010} it follows that the backstepping operator $\Pi:X\to X$ given by $\Pi(u_{0},u_{1})=(w_{0},w_{1})$, with $(w_{0},w_{1})$ defined in \eqref{transformed_initial}, has a linear and bounded inverse $\Pi^{-1}:X\to X$. We also mentioned Remark \ref{BS_operator}. Similarly to Remark \ref{BS_operator}, we also have that the inverse backstepping operator $\Pi^{-1}$ maps the state $(w,w_{t})$ of \eqref{wp6_dd} into the state $(u,u_{t})$ of \eqref{dd_1} with $U(t)=U_{1}(t)+U_{2}(t)$, where $U_{1}(t)$ is given by \eqref{control_1} and $U_{2}(t)$ is given by \eqref{m_energy_3}. Therefore, $(u,u_{t})\in W^{1,1}(0,\infty;X)$.\\

Let us note that $U_{1}(t)$ in \eqref{control_1} is written in terms of $(u,u_{t})$ but $U_{2}(t)$ in \eqref{m_energy_3} is written in terms of $(w,w_{t})$. However, we can rewrite $U_{2}(t)$ in \eqref{m_energy_3} in terms of $(u,u_{t})$ thanks to \eqref{transformed_initial}. Then, in order to present the final expression of the designed feedback law, for a regular enough function $f=f(t,x)$ let us introduce

\begin{equation}
\label{u1_dd}
U_{1}(t,f)=\frac{1}{h(L)}\left(\int_{0}^{L}k(L,y)f(t,y)\,dy+\int_{0}^{L}s(L,y)f_{t}(t,y)\,dy\right),
\end{equation}

\begin{equation}
\label{u2_dd}
U_{2}(t,f)=-\frac{1}{dh(L)}U_{3}(t,f)-P\,\text{sign}\left(U_{4}(t,f)\right),
\end{equation}

\begin{equation}
\label{u3_dd}
\begin{split}
U_{3}(t,f) =&~h(L)f_{t}(t,L)+s_{y}(L,L)f(t,L)-s(L,L)f_{x}(t,L)\\
 & -\int_{0}^{L}(2\lambda(y)s(L,y)+k(L,y))f_{t}(t,y)\,dy-\int_{0}^{L}(\beta(y)s(L,y)+s_{yy}(L,y))f(t,y)\,dy,
\end{split}
\end{equation}

\begin{equation}
\label{u4_dd}
\begin{split}
U_{4}(t,f) =&~h'(L)f(t,L)+h(L)f_{x}(t,L)-s(L,L)f_{t}(t,L)-k(L,L)f(t,L)\\
 & -\int_{0}^{L}s_{x}(L,y)f_{t}(t,y)\,dy-\int_{0}^{L}k_{x}(L,y)f(t,y)\,dy.
\end{split}
\end{equation}

\noindent Accordingly, in virtue of \eqref{transformed_initial} it follows from \eqref{u3_dd} that $w_{t}(t,L)=U_{3}(t,u)$, Also, in view of \eqref{BS} we see from \eqref{u4_dd} that $w_{x}(t,L)=U_{4}(t,u)$. Then, $U_{1}(t,u)=U_{1}(t)$ and $U_{2}(t,u)=U_{2}(t)$, where $U_{1}(t)$ and $U_{2}(t)$ are respectively given by \eqref{control_1} and \eqref{m_energy_3}.\\

Finally, the main result regarding the Dirichlet-Dirichlet case is the following one.

\begin{thm}
\label{dd_mr}
Let $\lambda\in C^{2}([0,L])$ and $\beta\in C([0,L])$. Let $d\in(0,\infty)$ be the desired decay rate. Let $k=k(x,y)$ and $s=s(x,y)$ be the kernel functions respectively obtained from \eqref{kernel_1} and \eqref{kernel_2}. Let us assume \textbf{(A1)} and \textbf{(A2)$\mathbf{_{DD}}$}. Let $(u_{0},u_{1})\in \Pi^{-1}(\mathcal{I}_{DD})=\{\Pi^{-1}(y_{0},y_{1})~/~(y_{0},y_{1})\in\mathcal{I}_{DD}\}$ be an initial condition, where $\mathcal{I}_{DD}=D(\mathcal{A})$ and $D(\mathcal{A})$ is given by \eqref{operator_dd}. Then, there exists a unique $u=u(t,x)$ with $(u,u_{t})\in W^{1,1}(0,\infty;X)$, where $X$ is given by \eqref{spaces}, such that

\begin{equation}
\label{dd_mr_closed}
\arraycolsep=0pt
\left\{\begin{array}{l}
u_{tt}=u_{xx}+2\lambda(x)u_{t}+\beta(x)u~\mbox{for almost every}~(t,x)\in(0,\infty)\times(0,L),\\[1ex]
u(t,0)=0~\mbox{for every}~t\in[0,\infty),\\[1ex]
u(t,L)\in p(t)+U_{1}(t,u)+U_{2}(t,u)~\mbox{for every}~t\in[0,\infty),\\[1ex]
u(0,x)=u_{0}(x)~\mbox{and}~u_{t}(0,x)=u_{1}(x)~\mbox{for every}~x\in(0,L).
\end{array}\right.
\end{equation}

\noindent Moreover, there exists a constant $C\in[1,\infty)$ such that

\begin{equation}
\label{dd_mr_decay}
\|u_{t}(t,\cdot)\|_{L^{2}(0,L)}^{2}+\|u_{x}(t,\cdot)\|_{L^{2}(0,L)}^{2}\leq Ce^{-2dt}\left(\|u_{1}\|_{L^{2}(0,L)}^{2}+\|u_{0}'\|_{L^{2}(0,L)}^{2}\right)\mbox{ for all }t\in[0,\infty).
\end{equation}
\end{thm}

\begin{rem}
Theorem \ref{dd_mr} can be applied to obtain \eqref{rapid} for \eqref{wave_dd}. To that end, it is necessary to take into account Remark \ref{simplification_dd} and make the necessary changes. Then, in view of Theorem \ref{kernel} (existence and uniqueness of regular enough kernel functions) we see that the additional assumption $\alpha\in C^{1}([0,L])$ is needed.
\end{rem}

\section{Dirichlet-Neumann case}

\hspace{0.75cm}In this section we prove Theorem \ref{dn_mr}, which is one of our main results and corresponds to the obtention of \eqref{rapid} for \eqref{wave_dn}. Most of the arguments in this section are mild variations to the ones presented previously for the Dirichlet-Dirichlet case.

\newpage
\mbox{}
\vspace{-1.25cm}

\begin{rem}
\label{simplification_dn}
We can transform \eqref{wave_dn} into another wave equation without first-order spatial derivative term by applying the transformation $v(t,x)=(e^{\frac{1}{2}\int_{0}^{x}\alpha(\tau)\,d\tau})u(t,x)$. Indeed, we get \eqref{dd_change} but with

\begin{equation*}
v_{x}(t,L)=(\alpha(L)/2)v(t,L)+(e^{\frac{1}{2}\int_{0}^{L}\alpha(\tau)\,d\tau})p(t)+(e^{\frac{1}{2}\int_{0}^{L}\alpha(\tau)\,d\tau})U(t)
\end{equation*}

\noindent instead of

\begin{equation*}
v(t,L)=(e^{\frac{1}{2}\int_{0}^{L}\alpha(\tau)\,d\tau})p(t)+(e^{\frac{1}{2}\int_{0}^{L}\alpha(\tau)\,d\tau})U(t).
\end{equation*}
\end{rem}

Thus, for $\gamma\in\mathbb{R}$, in this section we consider

\begin{equation}
\label{dn_1}
\left\{\begin{array}{rl}
u_{tt}=u_{xx}+2\lambda(x)u_{t}+\beta(x)u,&~(t,x)\in(0,\infty)\times(0,L),\\[1ex]
u(t,0)=0,~u_{x}(t,L)=\gamma u(t,L)+p(t)+U(t),&~t\in(0,\infty),\\[1ex]
u(0,x)=u_{0}(x),~u_{t}(0,x)=u_{1}(x),&~x\in(0,L).
\end{array}\right.
\end{equation}

\subsection{Feedback design}
\label{dn_fd}

\hspace{0.75cm}In this section we prove the feedback design part of Theorem \ref{dn_mr}. To that end, let us recall the required elements of the backstepping transformation constructed in \cite[Section 3]{sck2010}, which are: \eqref{BS_1}, \eqref{BS_2}, Theorem \ref{kernel} (existence and uniqueness of regular enough kernel functions) and \eqref{transformed_initial}.\\

As in the Dirichlet-Dirichlet case, in \eqref{dn_1} we consider $U(t)=U_{1}(t)+U_{2}(t)$. In order to choose $U_{1}(t)$ we first employ the backstepping transformation \eqref{BS} to get

\begin{equation*}
\begin{split}
w_{x}(t,x) =&~h'(x)u(t,x)+h(x)u_{x}(t,x)\\[0.5ex]
 & -k(x,x)u(t,x)-\int_{0}^{x}k_{x}(x,y)u(t,y)\,dy-s(x,x)u_{t}(t,x)-\int_{0}^{x}s_{x}(x,y)u_{t}(t,y)\,dy.
\end{split}
\end{equation*}

\noindent Then, in virtue of the computations done in \cite[Section 3.1]{sck2010} we choose

\begin{equation}
\label{control_2}
\begin{split}
U_{1}(t) =&-\gamma u(t,L)\\
 &+\frac{1}{h(L)}\left(-h'(L)u(t,L)+k(L,L)u(t,L)+\int_{0}^{L}k_{x}(L,y)u(t,y)\,dy\right.\\
 & \left.+s(L,L)u_{t}(t,L)+\int_{0}^{L}s_{x}(L,y)u_{t}(t,y)\,dy\right),
\end{split}
\end{equation}

\noindent to transform \eqref{dn_1} into the target system

\begin{equation}
\label{dn_2}
\left\{\begin{array}{rl}
w_{tt}=w_{xx}-2dw_{t}-d^{2}w,&~(t,x)\in(0,\infty)\times(0,L),\\[1ex]
w(t,0)=0,~w_{x}(t,L)=h(L)p(t)+h(L)U_{2}(t),&~t\in(0,\infty),\\[1ex]
w(0,x)=w_{0}(x),~w_{t}(0,x)=w_{1}(x),&~x\in(0,L),
\end{array}\right.
\end{equation}

\noindent where the transformed initial condition $(w_{0},w_{1})$ is given by \eqref{transformed_initial}.\\

Let us note that the formal computations performed in Section \ref{dd_fd} (Feedback design) for the Dirichlet-Dirichlet case are also valid for the Dirichlet-Neumann case. In particular, by keeping the definition of the modified energy \eqref{m_energy}, but now for a regular enough state $(w,w_{t})$ of the target system \eqref{dn_2}, we find that both \eqref{m_energy_1} and Lemma \ref{e_energy} are still valid.\\

As in the Dirichlet-Dirichlet case, we shall design $U_{2}(t)$ as a feedback law from \eqref{m_energy_1}. Let us assume \textbf{(A1)}. Then, as $w_{x}(t,L)=h(L)p(t)+h(L)U_{2}(t)$, it follows

\begin{align}
\notag
V'(t)+2dV(t) =&~U_{2}(t)h(L)\left(w_{t}(t,L)+dw(t,L)\right)+p(t)h(L)\left(w_{t}(t,L)+dw(t,L)\right)\\[1ex]
\label{dn_3}
 \leq&~U_{2}(t)h(L)\left(w_{t}(t,L)+dw(t,L)\right)+Ph(L)|w_{t}(t,L)+dw(t,L)|.
\end{align}

\noindent Then, with the choice of $U_{2}(t)$ so that

\begin{equation}
\label{dn_4}
U_{2}(t)=-P\,\text{sign}\left(w_{t}(t,L)+dw(t,L)\right),
\end{equation}

\noindent the right-hand side of \eqref{dn_3} becomes non-positive. Indeed, thanks to the property $\theta f=|f|$ for every $\theta\in\text{sign}(f)$, we have

\begin{align}
\notag
 &U_{2}(t)h(L)\left(w_{t}(t,L)+dw(t,L)\right)+Ph(L)|w_{t}(t,L)+dw(t,L)|\\[1ex]
\notag
 &=\theta h(L)\left(w_{t}(t,L)+dw(t,L)\right)+Ph(L)|w_{t}(t,L)+dw(t,L)|\\[1ex]
\notag
 &=-P h(L)|w_{t}(t,L)+dw(t,L)|+P h(L)|w_{t}(t,L)+dw(t,L)|\\[1ex]
\label{dn_5}
 &=0\mbox{ for every }\theta\in-P\,\text{sign}\left(w_{t}(t,L)+dw(t,L)\right).
\end{align}

\noindent Therefore, in virtue of \eqref{dn_3} and \eqref{dn_5} we conclude that the solutions of the corresponding closed-loop system, obtained by considering the target system \eqref{dn_2} together with \eqref{dn_4}, satisfy $V(t)\leq e^{-2dt}V(0)$ for all $t\in[0,\infty)$. The closed-loop system this time is the following one.

\begin{equation}
\label{dn_6}
\left\{\begin{array}{rl}
w_{tt}=w_{xx}-2dw_{t}-d^{2}w,&~(t,x)\in(0,\infty)\times(0,L),\\[1ex]
w(t,0)=0,&~t\in(0,\infty),\\[1ex]
w_{x}(t,L)\in h(L)p(t)-h(L)P\,\text{sign}\left(w_{t}(t,L)+dw(t,L)\right)
,&~t\in(0,\infty),\\[1ex]
w(0,x)=w_{0}(x),~w_{t}(0,x)=w_{1}(x),&~x\in(0,L).
\end{array}\right.
\end{equation}

To conclude, let us note that the arguments presented at the end of Section \ref{dd_fd} (Feedback design) for the Dirichlet-Dirichlet case regarding the backstepping operator $\Pi:X\to X$ and its inverse are also valid for the Dirichlet-Neumann case. Also:

\begin{rem}
\label{BS_operator_DN}
The backstepping operator $\Pi$ maps the state $(u,u_{t})$ of \eqref{dn_1} with \eqref{control_2} into the state $(w,w_{t})$ of the target system \eqref{dn_2}.
\end{rem}

Accordingly, \eqref{decay} is still valid, thus concluding \eqref{rapid} for \eqref{dn_1}. We have shown the feedback design part of Theorem \ref{dn_mr} and \eqref{dn_mr_decay} with $C=C_{1}C_{2}\sqrt{K}$, where $K$ is defined in \eqref{constant}.

\subsection{Well-posedness}
\label{dn_wp}

\hspace{0.75cm}In this section we apply the maximal monotone operator theory to prove the well-posedness of the closed-loop system \eqref{dn_6}, which is a differential inclusion, with the inclusion being located at the\l boundary condition. To that end, let us recall the Hilbert space \eqref{spaces} together with its inner product \eqref{inn-prod}.

\newpage

\noindent Let us introduce

\begin{equation}
\label{lift_dn}
y(t,x)=w(t,x)-\psi(t)\phi(x)~\mbox{with}~\psi(t)=h(L)p(t)~\mbox{and}~\phi(x)=\frac{x}{L}(x-L).
\end{equation}

\noindent Since $\phi(0)=0$, $\phi(L)=0$ and $\phi'(L)=1$ we get that $y(t,0)=w(t,0)$, $y(t,L)=w(t,L)$, $y_{t}(t,L)=w_{t}(t,L)$ and $y_{x}(t,L)=w_{x}(t,L)-h(L)p(t)$. Let us assume \textbf{(A2)$\mathbf{_{DN}}$}, thus obtaining that $y(0,x)=w(0,x)$ and $y_{t}(0,x)=w_{t}(0,x)$. Then, it follows from \eqref{dn_6} that $y=y(t,x)$ satisfies

\begin{equation}
\label{dn_7}
\left\{\begin{array}{rl}
y_{tt}=y_{xx}-2dy_{t}-d^{2}y+f,&~(t,x)\in(0,\infty)\times(0,L),\\[1ex]
y(t,0)=0,&~t\in(0,\infty),\\[1ex]
y_{x}(t,L)\in -h(L)P\,\text{sign}\left(y_{t}(t,L)+dy(t,L)\right),&~t\in(0,\infty),\\[1ex]
y(0,x)=w_{0}(x),~y_{t}(0,x)=w_{1}(x),&~x\in(0,L),
\end{array}\right.
\end{equation}

\noindent where $f(x,t) = (-\partial_{tt} +\partial_{xx} - 2d \partial_t - d^2)(\psi(t)\phi(x))$. Finally, let us consider the operator

\begin{equation}
\label{operator_dn}
\left\{\begin{array}{l}
\mathcal{A}:D(\mathcal{A})\subset X\to X,\\[2ex]
\mathcal{A}(q,l)=(-l,-q''+2dl+d^{2}q),\\[2ex]
D(\mathcal{A})=\{(q,l)\in X~/~\mathcal{A}(q,l)\in X,~q'(L)\in-h(L)P\,\text{sign}(l(L)+dq(L))\}.
\end{array}\right.
\end{equation}

\begin{rem}
As in Remark \ref{ope_rem_dd} we have that $D(\mathcal{A})$ is a non-empty set. Also, $\mathcal{A}$ is not a linear operator and $D(\mathcal{A})$ is dense in $X$. Finally,\\

$D(\mathcal{A})=\{(q,l)\in(H^{2}(0,L)\cap H)\times H~/~q'(L)\in-h(L)P\,\text{sign}(l(L)+dq(L))\}$.
\end{rem}

Thanks to \eqref{operator_dn} we can write \eqref{dn_7} in operator form as follows:

\begin{equation}
\label{operator_dn_2}
\left\{\begin{array}{rl}
\dfrac{d}{dt}(y,y_{t})+\mathcal{A}(y,y_{t})=(0,f),&t\in[0,\infty),\\[2ex]
(y,y_{t})(0)=(w_{0},w_{1}).
\end{array}\right.
\end{equation}

We proceed to prove that \eqref{operator_dn_2} is well-posed by applying the maximal monotone operator theory. In that direction we have the following two results: the first one (Proposition \ref{wp1_dn}) states that the operator $\mathcal{A}$ is monotone, while the second one (Proposition \ref{wp2_dn}) states that the operator $I+\mathcal{A}$ is surjective.

\begin{prop}
\label{wp1_dn}
The operator defined by \eqref{operator_dn} is monotone.
\end{prop}

\noindent\textbf{Proof~~}The proof is similar to that of Proposition \ref{wp1_dd} (Dirichlet-Dirichlet case). Let us employ the notation of that result. This time we also get that $\langle \mathcal{A} (q_1,l_1) - \mathcal{A}(q_2,l_2), (q_1,l_1)-(q_2,l_2) \rangle_{X}=I_{1}+I_{2}$ and conclude $I_{1}\geq0$ by the same means.\\

Let us handle $I_{2}=-l(L)q'(L)-dq(L)q'(L)$. As $q_{j}'(L)\in-h(L)P\,\text{sign}(l_{j}(L)+dq_{j}(L))$ for $j\in\{1,2\}$, there exists $\theta_{j}\in\text{sign}(l_{j}(L)+dq_{j}(L))$ such that $q_{j}'(L)=-h(L)P\theta_{j}$. Therefore,

\begin{align*}
 &I_{2}=-l(L)q'(L)-dq(L)q'(L)\\[1ex]
 &=-(q_{1}'(L)-q_{2}'(L))\left((l_{1}(L)+dq_{1}(L))-(l_{2}(L)+dq_{2}(L))\right)\\[1ex]
 &=h(L)P(\theta_{1}-\theta_{2})\left((l_{1}(L)+dq_{1}(L))-(l_{2}(L)+dq_{2}(L))\right).
\end{align*}

\noindent We conclude that $I_{2}\geq0$ because the sign multivalued operator \eqref{sign} is monotone. Accordingly, the operator $\mathcal{A}$ is monotone as $I_{1}+I_{2}\geq0$.\hfill$\blacksquare$

\begin{prop}
\label{wp2_dn}
The operator defined by \eqref{operator_dn} satisfies $R(I+\mathcal{A})=X$.
\end{prop}

\noindent\textbf{Proof~~}The proof is similar to that of Proposition \ref{wp2_dd} (Dirichlet-Dirichlet case). Given $(m,n)\in X$ we need to prove the existence of $(q,l)\in D(\mathcal{A})$, domain being defined in \eqref{operator_dn}, such that $(q,l)+\mathcal{A}(q,l)=(m,n)$. To that end, let us introduce the problem

\begin{equation}
\label{wp2_dn_1}
\left\{\begin{array}{l}
-q''+q(1+2d+d^{2})=n+m(1+2d),\\[1ex]
q(0)=0,~q'(L)\in-h(L)P\,\text{sign}((1+d)q(L)-m(L)).
\end{array}\right.
\end{equation}

\noindent By proving that \eqref{wp2_dn_1} possesses a unique solution $q\in H^{2}(0,L)$ and by putting $l=q-m$, then we get that $(q,l)\in D(\mathcal{A})$ and that $(q,l)+\mathcal{A}(q,l)=(m,n)$. As in the proof of Proposition \ref{wp2_dd} (Dirichlet-Dirichlet case), we analyze an optimization problem. Thus, let us introduce the functional $J:H\to\mathbb{R}$ by

\begin{align}
\notag
J(q)= &~\frac{1}{2}\int_{0}^{L}\left[(q')^{2}+q^{2}(1+2d+d^{2})-2q\left(n+m(1+2d)\right)\right]dx\\[1ex]
\label{functional_dn}
&+\frac{1}{1+d}h(L)P\,\varphi_{\sigma}\left((1+d)q(L)-m(L)\right),
\end{align}

\noindent where $\varphi_{\sigma}:\mathbb{R}\to\mathbb{R}$ is the Moreau regularization of the convex and continuous function $\varphi:\mathbb{R}\to\mathbb{R}$ given by $\varphi(x)=|x|$.\\

As in the proof of Proposition \ref{wp2_dd} (Dirichlet-Dirichlet case), we can argue the existence of a minimizer $q\in H$ for the functional $J$. Also, the G\^ateaux derivative of the functional $J$ at such minimizer must vanish for each direction $\rho\in H$. Accordingly, as $\varphi_{\sigma}'(x)=\alpha_{\sigma}(x)$ we obtain

\begin{align}
\notag
J'(q;\rho)= &~\int_{0}^{L}\left[q'\rho'+q\rho(1+2d+d^{2})-\rho\left(n+m(1+2d)\right)\right]dx\\[1ex]
\label{derivative_dn}
&+h(L)P\,\alpha_{\sigma}\left((1+d)q(L)-m(L)\right)\rho(L)=0~~\forall \rho\in H.
\end{align}

\noindent As $C_{0}^{\infty}(0,L)\subset H$, we can argue that $q\in H^{2}(0,L)$ and then that $-q''+q(1+2d+d^{2})=n+m(1+2d)$ for almost every $x\in(0,L)$. Moreover, after one integration by parts in \eqref{derivative_dn} we get

\begin{equation*}
\left[q'(L)+h(L)P\,\alpha_{\sigma}\left((1+d)q(L)-m(L)\right)\right]\rho(L)=0~~\forall \rho\in H.
\end{equation*}

Accordingly, so far we have shown:

\begin{lem}
\label{wp3_dn}
For any $\sigma\in(0,\infty)$ there exists $q_{\sigma}\in H^{2}(0,L)$ such that $-q_{\sigma}''+q_{\sigma}(1+2d+d^{2})=n+m(1+2d)$ for almost every $x\in(0,L)$, $q_{\sigma}(0)=0$ and $q_{\sigma}'(L)=-h(L)P\,\alpha_{\sigma}\left((1+d)q_{\sigma}(L)-m(L)\right)$.
\end{lem}

As in the proof of Proposition \ref{wp2_dd} (Dirichlet-Dirichlet case), let $q_{\sigma}$ be given by Lemma \ref{wp3_dn}. We proceed to obtain upper bounds, independent of $\sigma$, for $\|q_{\sigma}\|_{H^{2}(0,L)}$ and $|\alpha_{\sigma}\left((1+d)q_{\sigma}(L)-m(L)\right)|$.\\

$\bullet$ \underline{Upper bound for $\|q_{\sigma}\|_{H^{2}(0,L)}$ independent of $\sigma$}\\

By Lemma \ref{wp3_dn} we have

\begin{equation*}
\int_{0}^{L}\left[-q_{\sigma}''+q_{\sigma}(1+2d+d^{2})\right]\left[(1+d)q_{\sigma}-m\right]dx=\int_{0}^{L}\left[n+m(1+2d)\right]\left[(1+d)q_{\sigma}-m\right]dx.
\end{equation*}

\noindent Let us recall that $m\in H$, and hence, $m\in H^{1}(0,L)$ with $m(0)=0$. Then, after one integration by parts and applying Cauchy inequality we get

\begin{align}
\notag
&(1+d)\int_{0}^{L}|q_{\sigma}'|^{2}dx-q_{\sigma}'(L)\left[(1+d)q_{\sigma}(L)-m(L)\right]-\int_{0}^{L}q_{\sigma}'m'\,dx\\[1ex]
\notag
&+(1+d)(1+d)^{2}\int_{0}^{L}|q_{\sigma}|^{2}dx-(1+d)^{2}\int_{0}^{L}q_{\sigma}m\,dx\\[1ex]
\label{uniform_dn}
&\leq\frac{1}{4}(1+d)^{2}\int_{0}^{L}|q_{\sigma}|^{2}dx+\|n+m(1+2d)\|_{L^{2}(0,L)}^{2}+\|n+m(1+2d)\|_{L^{2}(0,L)}\|m\|_{L^{2}(0,L)}.
\end{align}

\noindent We shall use \eqref{uniform_dn} later. Once again applying Cauchy inequality we arrive at

\begin{equation}
\label{uniform2_dn}
\int_{0}^{L}q_{\sigma}'m'\,dx\leq\frac{1}{2}\int_{0}^{L}|q_{\sigma}'|^{2}dx+\frac{1}{2}\|m'\|_{L^{2}(0,L)}^{2},
\end{equation}

\begin{equation}
\label{uniform3_dn}
(1+d)^{2}\int_{0}^{L}q_{\sigma}m\,dx\leq\frac{1}{2}(1+d)^{2}\int_{0}^{L}|q_{\sigma}|^{2}dx+\frac{1}{2}(1+d)^{2}\|m\|_{L^{2}(0,L)}^{2}.
\end{equation}

\noindent In order to treat the term $-q_{\sigma}'(L)\left[(1+d)q_{\sigma}(L)-m(L)\right]$ we use that $\alpha_{\sigma}(0)=0$. Then, by Lemma \ref{wp3_dn} and as $\alpha_{\sigma}$ is a monotone operator, we infer

\begin{align}
\notag
&-q_{\sigma}'(L)\left[(1+d)q_{\sigma}(L)-m(L)\right]\\[1ex]
\notag
&=h(L)P\left[\alpha_{\sigma}((1+d)q_{\sigma}(L)-m(L))\right]\left[(1+d)q_{\sigma}(L)-m(L)\right]\\[1ex]
\label{uniform4_dn}
&=h(L)P\left[\alpha_{\sigma}((1+d)q_{\sigma}(L)-m(L))-\alpha_{\sigma}(0)\right]\left[((1+d)q_{\sigma}(L)-m(L))-0\right]\geq0.
\end{align}

\noindent Accordingly, from the combination of \eqref{uniform_dn}, \eqref{uniform2_dn}, \eqref{uniform3_dn} and \eqref{uniform4_dn} we can get an upper bound for $\|q_{\sigma}\|_{H^{1}(0,L)}$, independent of $\sigma$, in terms of $\|m\|_{H^{1}(0,L)}$ and $\|n\|_{L^{2}(0,L)}$. Furthermore, thanks to Lemma \ref{wp3_dn} we can get the same for $\|q_{\sigma}\|_{H^{2}(0,L)}$. The conclusion is: there exists a positive constant $C_{1}=C_{1}(d)$, independent of $\sigma$, such that

\begin{equation}
\label{uniform5_dn}
\|q_{\sigma}\|_{H^{2}(0,L)}\leq C_{1}\|m\|_{H^{1}(0,L)}+C_{1}\|n\|_{L^{2}(0,L)}.
\end{equation}

$\bullet$ \underline{Upper bound for $|\alpha_{\sigma}\left((1+d)q_{\sigma}(L)-m(L)\right)|$ independent of $\sigma$}\\

There exists a positive constant $C_{2}=C_{2}(L)$ such that $\|q\|_{C^{1}([0,L])}\leq C_{2}\|q\|_{H^{2}(0,L)}$ for any $q\in H^{2}(0,L)$ because the injection of $H^{2}(0,L)$ into $C^{1}([0,L])$ is continuous. Then, by Lemma \ref{wp3_dn} and \eqref{uniform5_dn} we obtain

\begin{equation}
\label{uniform6_dn}
|\alpha_{\sigma}\left((1+d)q_{\sigma}(L)-m(L)\right)|\leq\frac{1}{h(L)P}|q_{\sigma}'(L)|\leq\frac{C_{2}C_{1}}{h(L)P}\|m\|_{H^{1}(0,L)}+\frac{C_{2}C_{1}}{h(L)P}\|n\|_{L^{2}(0,L)}.
\end{equation}

Therefore, in view of \eqref{uniform5_dn} and \eqref{uniform6_dn} we see that the sequences

\begin{equation*}
\left\{q_{\sigma}\right\}_{\sigma>0}\subset H^{2}(0,L)\mbox{ and }\left\{\alpha_{\sigma}\left((1+d)q_{\sigma}(L)-m(L)\right)\right\}_{\sigma>0}\subset\mathbb{R}
\end{equation*}

\noindent are bounded. Then, there exist $(q,b)\in H^{2}(0,L)\times\mathbb{R}$ and subsequences, which we denote by the same symbols, such that $q_{\sigma}\rightharpoonup q$ in $H^{2}(0,L)$ (weak convergence) and $\alpha_{\sigma}\left((1+d)q_{\sigma}(L)-m(L)\right)\to b$ in $\mathbb{R}$ as $\sigma\to0^{+}$. Moreover, we infer that $q_{\sigma}\to q$ in $C^{1}([0,L])$ as $\sigma\to0^{+}$ since the injection of $H^{2}(0,L)$ into $C^{1}([0,L])$ is compact, implying that $q_{\sigma}(x)\to q(x)$ and $q_{\sigma}'(x)\to q'(x)$ in $\mathbb{R}$ for all $x\in[0,L]$ as $\sigma\to0^{+}$.\\

Accordingly, Lemma \ref{wp3_dn} and the previous arguments yield:

\begin{lem}
\label{wp4_dn}
There exists $q\in H^{2}(0,L)$ such that $-q''+q(1+2d+d^{2})=n+m(1+2d)$ for almost every $x\in(0,L)$, $q(0)=0$ and $q'(L)=-h(L)Pb$.
\end{lem}

In view of Lemma \ref{wp4_dn} and \eqref{wp2_dn_1} we see that in order to complete the proof of $R(I+\mathcal{A})=X$ we just need to prove that $b\in\alpha\left((1+d)q(L)-m(L)\right)=\text{sign}\left((1+d)q(L)-m(L)\right)$. The required arguments can be found at the end of the proof of \cite[Proposition 3.1]{gh2023} or at the end of the proof of \cite[Proposition 3.2]{gh2025}. The proof of Proposition \ref{wp2_dn} is complete.\hfill$\blacksquare$\\

Let us return to \eqref{operator_dn_2}. We argue as we did at the end of Section \ref{dd_wp} (Well-posedness) for the Dirichlet-Dirichlet case. We have that the operator $\mathcal{A}$ defined in \eqref{operator_dn} is maximal monotone. Let us recall \textbf{(A2)$\mathbf{_{DN}}$} and let us assume $(w_{0},w_{1})\in D(\mathcal{A})$. Then, we conclude the existence of a unique $(y,y_{t})\in W^{1,1}(0,\infty;X)$ such that

\begin{equation}
\label{wp5_dn}
\left\{\begin{array}{l}
\dfrac{d}{dt}(y,y_{t})+\mathcal{A}(y,y_{t})=(0,f)\mbox{ for almost every }t>0,\\[2ex]
(y,y_{t})\in D(\mathcal{A})\mbox{ for every }t\geq0,\\[2ex]
(y,y_{t})(0)=(w_{0},w_{1}).
\end{array}\right.
\end{equation}

\noindent With \eqref{dn_6} (the closed-loop system), \eqref{lift_dn} and \eqref{dn_7} in mind, due to \eqref{wp5_dn} we conclude that $(w,w_{t})\in W^{1,1}(0,\infty;X)$ uniquely solves

\begin{equation}
\label{wp6_dn}
\arraycolsep=0pt
\left\{\begin{array}{l}
w_{tt}=w_{xx}-2dw_{t}-d^{2}w ~\mbox{for almost every}~(t,x)\in(0,\infty)\times(0,L),\\[1ex]
w(t,0)=0~\mbox{for every}~t\in[0,\infty),\\[1ex]
w_{x}(t,L)\in h(L)p(t)-h(L)P\,\text{sign}\left(w_{t}(t,L)+dw(t,L)\right)~\mbox{for every}~t\in[0,\infty),\\[1ex]
w(0,x)=w_{0}(x)~\mbox{and}~w_{t}(0,x)=w_{1}(x)~\mbox{for every}~x\in(0,L).
\end{array}\right.
\end{equation}

We have shown the well-posedness part of Theorem \ref{dn_mr}.

\subsection{Conclusion}
\label{dn_c}

\hspace{0.75cm}Here we gather the conclusions of both the feedback design (see Section \ref{dn_fd}) and the well-posedness (see Section \ref{dn_wp}) parts to establish our main result regarding the Dirichlet-Neumann case. The arguments are similar to that of Section \ref{dd_c} (Dirichlet-Dirichlet case).\\

From \cite[Section 3.3]{sck2010} it follows that the backstepping operator $\Pi:X\to X$ given by $\Pi(u_{0},u_{1})=(w_{0},w_{1})$, with $(w_{0},w_{1})$ defined in \eqref{transformed_initial}, has a linear and bounded inverse $\Pi^{-1}:X\to X$. We also mentioned Remark \ref{BS_operator_DN}. Similarly to Remark \ref{BS_operator_DN}, we also have that the inverse backstepping operator $\Pi^{-1}$ maps the state $(w,w_{t})$ of \eqref{wp6_dn} into the state $(u,u_{t})$ of \eqref{dn_1} with $U(t)=U_{1}(t)+U_{2}(t)$, where $U_{1}(t)$ is given by \eqref{control_2} and $U_{2}(t)$ is given by \eqref{dn_4}. Therefore, $(u,u_{t})\in W^{1,1}(0,\infty;X)$.\\

Let us note that $U_{1}(t)$ in \eqref{control_2} is written in terms of $(u,u_{t})$ but $U_{2}(t)$ in \eqref{dn_4} is written in terms of $(w,w_{t})$. However, we can rewrite $U_{2}(t)$ in \eqref{dn_4} in terms of $(u,u_{t})$ thanks to \eqref{transformed_initial}. Then, in order to present the final expression of the designed feedback law, for a regular enough function $f=f(t,x)$ let us introduce

\begin{equation}
\label{u1_dn}
\begin{split}
U_{1}(t,f) =&-\gamma f(t,L)\\
&+\frac{1}{h(L)}\left(-h'(L)f(t,L)+k(L,L)f(t,L)+\int_{0}^{L}k_{x}(L,y)f(t,y)\,dy\right.\\
 & \left.+s(L,L)f_{t}(t,L)+\int_{0}^{L}s_{x}(L,y)f_{t}(t,y)\,dy\right),
\end{split}
\end{equation}

\begin{equation}
\label{u2_dn}
U_{2}(t,f)=-P\,\text{sign}\left(U_{3}(t,f)+dU_{4}(t,f)\right),
\end{equation}

\begin{equation}
\label{u3_dn}
\begin{split}
U_{3}(t,f) =&~h(L)f_{t}(t,L)+s_{y}(L,L)f(t,L)-s(L,L)f_{x}(t,L)\\
 & -\int_{0}^{L}(2\lambda(y)s(L,y)+k(L,y))f_{t}(t,y)\,dy-\int_{0}^{L}(\beta(y)s(L,y)+s_{yy}(L,y))f(t,y)\,dy,
\end{split}
\end{equation}

\begin{equation}
\label{u4_dn}
U_{4}(t,f)=h(L)f(t,L)-\int_{0}^{L}k(L,y)f(t,y)\,dy-\int_{0}^{L}s(L,y)f_{t}(t,y)\,dy.
\end{equation}

\noindent Accordingly, in virtue of \eqref{transformed_initial} it follows from \eqref{u3_dn} that $w_{t}(t,L)=U_{3}(t,u)$, Also, in view of \eqref{BS} we see from \eqref{u4_dn} that $w(t,L)=U_{4}(t,u)$. Then, $U_{1}(t,u)=U_{1}(t)$ and $U_{2}(t,u)=U_{2}(t)$, where $U_{1}(t)$ and $U_{2}(t)$ are respectively given by \eqref{control_2} and \eqref{dn_4}.\\

Finally, the main result regarding the Dirichlet-Neumann case is the following one.

\begin{thm}
\label{dn_mr}
Let $\lambda\in C^{2}([0,L])$, $\beta\in C([0,L])$ and $\gamma\in\mathbb{R}$. Let $d\in(0,\infty)$ be the desired decay rate. Let $k=k(x,y)$ and $s=s(x,y)$ be the kernel functions respectively obtained from \eqref{kernel_1} and \eqref{kernel_2}. Let us assume \textbf{(A1)} and \textbf{(A2)$\mathbf{_{DN}}$}. Let $(u_{0},u_{1})\in \Pi^{-1}(\mathcal{I}_{DN})=\{\Pi^{-1}(y_{0},y_{1})~/~(y_{0},y_{1})\in\mathcal{I}_{DN}\}$ be an initial condition, where $\mathcal{I}_{DN}=D(\mathcal{A})$ and $D(\mathcal{A})$ is given by \eqref{operator_dn}. Then, there exists a unique $u=u(t,x)$ with $(u,u_{t})\in W^{1,1}(0,\infty;X)$, where $X$ is given by \eqref{spaces}, such that

\begin{equation}
\label{dn_mr_closed}
\arraycolsep=0pt
\left\{\begin{array}{l}
u_{tt}=u_{xx}+2\lambda(x)u_{t}+\beta(x)u~\mbox{for almost every}~(t,x)\in(0,\infty)\times(0,L),\\[1ex]
u(t,0)=0~\mbox{for every}~t\in[0,\infty),\\[1ex]
u_{x}(t,L)\in\gamma u(t,L)+p(t)+U_{1}(t,u)+U_{2}(t,u)~\mbox{for every}~t\in[0,\infty),\\[1ex]
u(0,x)=u_{0}(x)~\mbox{and}~u_{t}(0,x)=u_{1}(x)~\mbox{for every}~x\in(0,L).
\end{array}\right.
\end{equation}

\noindent Moreover, there exists a constant $C\in[1,\infty)$ such that

\begin{equation}
\label{dn_mr_decay}
\|u_{t}(t,\cdot)\|_{L^{2}(0,L)}^{2}+\|u_{x}(t,\cdot)\|_{L^{2}(0,L)}^{2}\leq Ce^{-2dt}\left(\|u_{1}\|_{L^{2}(0,L)}^{2}+\|u_{0}'\|_{L^{2}(0,L)}^{2}\right)\mbox{ for all }t\in[0,\infty).
\end{equation}
\end{thm}

\begin{rem}
Theorem \ref{dn_mr} can be applied to obtain \eqref{rapid} for \eqref{wave_dn}. To that end, it is necessary to take into account Remark \ref{simplification_dn}, put $\gamma=\alpha(L)/2$ and make the necessary changes. Then, in view of Theorem \ref{kernel} (existence and uniqueness of regular enough kernel functions) we see that the additional assumption $\alpha\in C^{1}([0,L])$ is needed.
\end{rem}

\section{Concluding remarks}
\label{cr}

\hspace{0.75cm}In this paper we have shown the rapid stabilization \eqref{rapid} for the unstable wave equations \eqref{wave_dd} and \eqref{wave_dn}, in which there is an unknown disturbance acting at the boundary condition. We have employed the backstepping method, Lyapunov techniques and the sign multivalued operator \eqref{sign} for the design of suitable feedback laws. The assumptions made on the unknown disturbance, namely \textbf{(A1)} and \textbf{(A2)}, are the standard ones that can be found in the literature. The well-posedness of the corresponding closed-loop system has been shown with the maximal monotone operator theory. The main limitation of the approach followed in this paper is that it is needed a backstepping transformation to map the original plant into a target system so that Lyapunov techniques can be applied.

\section{Acknowledgments}

\hspace{0.75cm}Patricio Guzm\'an and Agust\'in Huerta received financial support from FONDECYT 11240290. Hugo Parada was partially supported by the Agence Nationale de la Recherche through the Labex CIMI (ANR-11-LABX-0040). He is currently supported by the QuBiCCS project (ANR-24-CE40-3008).


\bibliographystyle{plain}
\bibliography{references}


\end{document}